\documentclass{article}

    \usepackage[preprint]{neurips_2021}

\usepackage[utf8]{inputenc} 
\usepackage[T1]{fontenc}    
\usepackage{url}            
\usepackage{booktabs}       
\usepackage{amsfonts}       
\usepackage{nicefrac}       
\usepackage{microtype}      
\usepackage{xcolor}         
\usepackage{amsmath}
\usepackage{capt-of}
\usepackage{graphicx}

\title{SCIMAT: Science and Mathematics Dataset\footnote{Extended version is submitted to archival journal.}}

\author{%
  Neeraj Kollepara \\
    CSTAR, IIIT, 
  Hyderabad \\
  \texttt{neeraj.kollepara@students.iiit.ac.in} \\
  \And
  Snehith Kumar Chatakonda \\
    CSTAR, IIIT,
  Hyderabad \\
  \texttt{snehith.kumar@students.iiit.ac.in} \\
   \And
   Pawan Kumar \\
  CSTAR, IIIT,
  Hyderabad \\
  \texttt{pawan.kumar@iiit.ac.in} \\
}

\begin{document}

\maketitle

\begin{abstract}
  In this work, we announce a comprehensive well curated and opensource dataset with millions of samples for pre-college and college level problems in mathematics and science. A preliminary set of results using transformer architecture with character to character encoding is shown. The dataset identifies some challenging problem and invites research on better architecture search. 
\end{abstract}

\section{Introduction and Previous Work}
Datasets play an important role in driving research in supervised machine learning research. Some prominant examples being MNIST \cite{lecun-mnisthandwrittendigit-2010} for hand written digit classification, CIFAR10 \cite{cifar10} and IMAGENET \cite{imagenet} for image classification and generative models, etc. For solving math word problems, semantic rules and various models have been proposed in NLP community since 1963 starting from \cite{bobrow}, \cite{briars1984}, \cite{feigenbaum1963computers}. Some word problems in \cite{fletcher1985} are of the form: {\tt  Lucy has two dimes. Sarah has six dimes. How many dimes do they have altogether?} or {\tt Dan has six books. Jill has two books. How many
books does Dan have more than Jill?} This paper uses Kintsch and Greeno's (1985) theory of comprehension and solution for arithmetic word problems above. These papers used classical approaches with semantic rules. Recently, machine learning models have been used for which large labelled dataset is essential. Hence, there is a dire need of large question-answer dataset for mathematics and science problems; such dataset can have impact on online education, intelligent tutoring and automated grading. For intelligent tutoring, not just the answers, but the step by step hint can be provided; this is explored in \cite{kang2016}. However, tutoring requires some knowledge graph representation. Although, this was shown for simple algebraic and geometric mathematics problems, it remains a challenging task for more advanced problems. No wonder tutoring is a complex task as nicely pointed out in detail in \cite{kenneth2013}. Given that intelligent tutoring is one of the most challenging task, the datasets and innovative architectures would play a critical role to succeed in this endevour. Recently, question answer dataset\footnote{\url{https://github.com/deepmind/mathematics\_dataset}} for mathematics was proposed in \cite{AnalysingMR} and for word problem sample dataset was proposed in \cite{SMWPSDL1}, and a comparison of results for character to character encoding for transformer and for LSTM is shown. This dataset has selected problems in mathematics for math exams for British 16 year old school children. Some sample questions are: ${\tt Factorise~ x^2 + 7x}$ or {\tt Three letters picked without replacement from qqqkkklkqkkk. Give prob of sequence qql.}  In \cite{clark2018think}, a set of 7787 multiple choice questions in high school science questions is proposed as ARC (AI2 Reasoning challenge, 2018). A sample question from this dataset is: {\tt Which property of a mineral can be determined just by looking at it? (A) luster [correct] (B) mass (C) weight (D) hardness}. Moreover, with the ARC challenge a large corpus of 14 million science sentences relevant to the question-answer set is also proposed. A sample sentence from the corpus is: {\tt Random motion of the air molecules and turbulence provide upward forces that may counteract the downward force of gravity.} Such a corpus allows language understanding and questions with linguistic variations. We remark that any other corpus can be used for training the given architecture for linguistic understandings, which is further trained on the given datasets. For the ARC challenge, several baseline neural models were proposed. There are datasets for logical reasoning and English comprehension. For example, in \cite{weston2015aicomplete}, logical reasoning question answer dataset is proposed. The reasoning is considered to be of various types such as problems involving single supporting fact, two supporting fact, counting, path finding, size reasoning, etc. A sample question for path finding is: {\tt The kitchen is north of the hallway. John is hungry. The bathroom is west of the bedroom. John goes to the kitchen.
The den is east of the hallway. John grabbed the apple there. The office is south of the bedroom. Daniel is hungry. How do you go from den to kitchen? How do you go from office to bathroom?}. The last two sentences are questions with answers {\tt west, north} and {\tt north, west} respectively. This dataset is part of bAbI project\footnote{\url{https://github.com/facebookarchive/bAbI-tasks}} of facebook research. For algebra word problems, a dataset\footnote{\url{http://groups.csail.mit.edu/rbg/code/wordprobs/}} and code is proposed in \cite{kushman-etal-2014-learning}. Most of these word problems correspond to solving system of linear equations, their method derives these equations, then solves it. A sample question answer in this dataset taken from \cite{kushman-etal-2014-learning} is: {\tt An amusement park sells 2 kinds of tickets. Tickets for children cost \$1.50. Adult tickets cost \$4. On a certain day, 278 people entered the park. On that same day the admission fees
collected totaled \$792. How many children were admitted on that day? How many adults
were admitted?} with soutions {\tt x = 128, y = 150}. Continuing along these lines in \cite{wang-etal-2017-deep}, they propose to translate the math word problem to equation using recurrent neural network (RNN) without doing any complex feature extractions. To the best of our knowledge, a comprehensive {\bf opensource} dataset for mathematics and science for pre-college and college level have been missing. To this end, in the following, we announce a new large dataset named SCIMAT, and we show preliminary results and comparisons.  

\section{SCIMAT: Large Science and Mathematics dataset}
We announce a large dataset\footnote{\url{https://github.com/misterpawan/scimat2}} of hundreds of millions of question-answer for mathematics and science for pre-college and college level, which typically is taught to 15-19 age group around the world. The list of topics covered are: Acids And Bases, Atomic Structure, Stoichiometry, Thermodynamics, Units And Dimensions, Kinematics, Laws of Motion, Work Power Energy, Rotatory Motion, Gravitation, Electricity, Moving Charges and Magnetism, Electro Magnetic Induction, Alternating Current, Electro Magnetic Waves,  Ray Optics and Optical Instruments, Wave Optics, Dual Nature of Matter,  Mechanical Properties of Solids and Liquids,  Thermal Properties of Matter, Kinetic theory of Gases, Sound, Waves And Oscillations, SemiConductors, Communication Systems, etc. 
Each topic contains several subtopics, where each subtopics has hundreds of thousands of question answer dataset. 
\subsection{Sample Questions in Science}
\small
\begin{enumerate}
    \item {\bf Question:} 33 mL of a solution of HNO3 is found to be completely neutralised by 45 mL of a given solution of NaOH. If we take 12 mL of the same solution of HNO3, the amount of NaOH solution (the same solution as before) required to neutralise it will be.
    {\bf Answer:} 16.36 ml
 
    \item {\bf Question:} If a diatomic gas of 1 moles at 68 atm and volume 68 lit is adiabatically changed to volume 188 lit, then what will be the pressure.
    {\bf Answer :} 16.4atm
 
    \item {\bf Question:} A body is dropped from a height of 9578 m with an initial velocity of 42 m/s. With what velocity will it strike the ground ?
    {\bf Answer:} 435.3 m/s

    \item {\bf Question:} A 9062 N force is applied on a body of mass 980 kg placed on a smooth surface, then what is the resulting acceleration obtained ?
    {\bf Answer:} 9.2 m/s2

    \item {\bf Question:} The volume of 549 g of a substance is 116 cm3. If the density of liquid in which substance is placed is 4 g/cm3, will the substance float or sink ?
    {\bf Answer:} sink
    
    \item {\bf Question:} If a 822 V battery is connected across an unknown resistor, there is 224 A in the circuit, find the value of resistance of the resistor ? 
    {\bf Answer:} 3.7 ohm
    
    \item{\bf Question:} A square coil of side 3 cm consists of 31 turns and carries a current of 5 A. The coil is suspended vertically and the normal to the plane of the coil makes an angle of 53 degress with the direction of a uniform horizontal magnetic field of magnitude 17 tesla. What is the magnitude of the torque experienced by the coil.
    {\bf Answer:} 1.9 newton-m
 
    \item {\bf Question:} A series LCR circuit is connected to a variable frequency 230 V source with L = 193 H, C = 72 muF, R = 176 ohm. Determine the rms potential drop across resistance?
    {\bf Answer:} 230 volt

    \item {\bf Question:} Suppose that the electric field amplitude of an electromagnetic wave is E0 = 1936 N/C and that its frequency is v = 1512 MHz. Find an expression for B?
    {\bf Answer:} 6.45e-06sin[3.17e+01x-9.50e+09t]
 
    \item {\bf Question:} During blood transfusion, the needle is inserted in a vein where the gauge pressure is 1720 Pa. If the blood container is placed at 177 mm above the earth level so that blood may just enter the vein, is it safe for the patient?.
    {\bf Answer:} yes, patient is safe
  
    \item {\bf Question:} A sound wave travels at a speed of 29980.8 m/s, if it's wavelength is 32 m, will the sound wave be audible ?
    {\bf Answer:} audible

    \item {\bf Question:} For an amplitude modulated wave, the maximum amplitude is found to be 18.62 V while the minimum amplitude is found to be 7.91 V. Determine the modulation index.
    {\bf Answer:} 0.4

\end{enumerate}
\normalsize

Similarly, for mathematics, we append datasets from calculus (differentiation and integration), linear algebra (rank, row reduced echelon form, determinant, trace, etc), set operations, statistics, number theory, probability, etc. Some sample questions in this dataset are the following:
\subsection{Sample Question in Mathematics}
\small
\begin{enumerate}
   
\item {\bf Question:} Differentiate 293 * x * (sin(x) + sec(x)) with respect to x
 \par
{\bf Answer:} 293 * x * (cos(x) + tan(x) * sec(x)) + 293 * sin(x) + 293 * sec(x) 
\par

\item {\bf Question:} Integrate cot(4*x \textsuperscript{$\wedge$}2) + sec(22*x\textsuperscript{$\wedge$}2) with respect to x
\par
{\bf Answer:} 8*x*( -cot(4 * x\textsuperscript{$\wedge$}2 ) \textsuperscript{$\wedge$}2 - 1) + 44*x* tan(22*x\textsuperscript{$\wedge$}2 ) * sec(22*x\textsuperscript{$\wedge$}2 ) \par

\item {\bf Question:} 2 * f ( x ) + 8 * Derivative ( f ( x ) , x ) + Derivative ( f ( x ) , ( x, 2 ) ) = 0
\par
{\bf Answer:} f ( x ) = ( C1 * exp ( x * ( 1 - sqrt ( 6 ) ) ) + C2 * exp ( x * ( 1 + sqrt ( 6 ) ) ) )

\item {\bf Question:} Calculate the Rank of Matrix  [ [2, 1, 3, 7] , [1, 0, 4, 2 ] , [ 3, 1, 7, 9 ] ] 
{\bf Answer:} 2 \par

 \item {\bf Question:} Calculate the Trace of Matrix  [ [ 13, 38, 61 ] , [ 29, 1, 39 ] , [ 92, 16, 45 ] ] 
 {\bf Answer:} 59 \par

\item  {\bf Question:} What is the union of  \{ 2, 6, 7, 8, 9 \}  with \{ 3, 7, 8 \} 
{\bf Answer: } \quad \{ 2, 3, 6, 7, 8, 9 \}

\item {\bf Question:} What is the median of the sequence  ( 20, 38, 4, 21, 31, 94, 55)
{\bf Answer: } 31

\item {\bf Question:} What is 2 (base 3) in base 8?
{\bf Answer:} 2

\item {\bf Question:} Expand (-s + s + 2*s**5)*(4 - 1 - 2) - 3*s**5 + 4*s**5 + 0*s**5 - 2*s**5 - s**5 + 5*s**5 + (3*s**2 - 4 + 4)*(5*s**3 - 5*s**3 - s**3).
{\bf Answer:} 2*s**5

\item {\bf Question:} Three letters picked without replacement from {a: 3, c: 1, b: 7, d: 3}. Give prob of sequence bdc.
{\bf Answer:} 1/104

\end{enumerate}

\section{Numerical Experiments}

\small 
\begin{table*}
\parbox{.45\linewidth}{
\centering
\begin{tabular}{ll}
\cmidrule(r){1-2}
{Type of problem}                                                 & 
\begin{tabular}[c]{@{}l@{}}{C2C }\\{Accuracy} \end{tabular}  \\
\midrule
Differentiation of sum   & 99\% \\ 
Differentiation of product  & 100\%  \\ 
Differentiation of composition  &  100\% \\ 
Integration of sum  & 100\% \\ 
Integration of product & 100\%  \\ 
Integration of composition  &   92.5\% \\
Addition of matrices   &  49\%  \\
Subtraction of matrices   &  74\% \\
Transpose of matrix   &  100\%  \\
Determinant of matrix    &  32\%  \\
Multiplication of matrices   &  32\%  \\
Trace of a matrix   &  100\%   \\
Product of matrix with a Scalar    &  100\%  \\
Row Reduced echelon   &  76\%    \\
Rank of a matrix  &  92.5\%    \\
Mean of a sequence   &  95\%\\
Variance of a sequence   &  39.5\%   \\
Median of a sequence   &  99\%  \\
Set Union   &  100\%  \\
Set intersection   &  97.5\%\\
Set difference    &  100\%  \\
Symmetric difference between sets   &  100\%  \\
\bottomrule
\end{tabular}
\caption{Comparison of our model trained on new datasets with Char2Char transformer. The Char2Char is denoted by C2C.}
\label{tab:new}
}
\hfill
\parbox{.45\linewidth}{

\begin{tabular}{ll}
\toprule
\cmidrule(r){1-2}
{Type of problem}                                                 & 
\begin{tabular}[c]{@{}l@{}}{C2C }\\{Accuracy} \end{tabular}  \\
\midrule
Neutralization  & 82.6\% \\
Adiabatic  & 76.3\% \\
Refrigrator & 82.6\% \\
Estimated value  & 61.8\% \\
Force, mass, acceleration  & 45.5\% \\ 
Momentun conservation  & 78.7\% \\
Kinetic energy  & 73.5\% \\ 
Balancing a metre stick  & 81.5\% \\
Gravitational field  & 94.2\% \\
Float or sink?  &  98\% \\ 
Ohms Law  & 89.0\% \\
Torque due to magnetic field  & 84.5\%\\
LCR circuit  & 91.3\% \\
Mirror formula for concave & 79.6\% \\
Is the sound audible?  & 76\% \\
Sound wave propogation  &   31.5\% \\
modulation index  & 89.6\% \\
Force between wires & 75.2\% \\
Conservation of momentum  & 14.5\%  \\ 
Potential energy & 63\%  \\ 
Work, mass, velocity & 10\% \\ 
\bottomrule
\end{tabular}
\caption{Accuracy for science datasets with Char2Char transformer. The Char2Char is denoted by C2C. See dataset with folder names.}
\label{tab:science}
}
\end{table*}
\normalsize

The  code  for training and testing is  written  in  Python  and  {\tt PyTorch} framework  is  used.  The  models  are  trained on dual {\tt Intel Xeon E5-2640 v4} processors, providing 40 virtual cores per node, 128 GB of 2400MT/s DDR4 ECC RAM and four {\tt Nvidia GeForce GTX 1080} TiGPUs, providing 14336 CUDA cores.
We use the standard transformer described in \cite{vaswani2017} with our own specifications as follows. We use an encoder which is composed of stack of $N = 4$ identical layers. The embedding size (dmodel) = 128, attention heads ($h$) = 8. The inner layer size of feed forward network used in each layers of encoder stack (dff) = 512. We minimize the sum of log probabilities of the correct tokens via the Adam optimizer with adaptive learning rate. The model was trained for 100 epochs. For floating point answers, accuracy for two digits after decimal place was matched. In Table \ref{tab:new}, \ref{tab:science}, we find that there are datasets where it is challenging to obtain high accuracy, and robust architecture or encoding is required. Since lately, many other variants of transformers were proposed, in Table \ref{compare}, we compare various different transformer with word-to-word and char-to-char encoding. In general, we found that char2char gives best accuracy. 

\small
\begin{table}[t]
\centering
\setlength{\tabcolsep}{1.249em}
\begin{tabular}{@{}lccc@{}}
\toprule
\multicolumn{1}{c}{Type of problem} & {\begin{tabular}[c]{@{}c@{}}C2C Trans.\end{tabular}} & {\begin{tabular}[c]{@{}c@{}}W2W Trans.\end{tabular}} & {\begin{tabular}[c]{@{}c@{}}W2W Perfor.\end{tabular}} \\ 
\midrule
pH &  \textbf{99.8\%} & 97.3\% & 84.7\% \\
Compare number of atoms  & \textbf{97.5\%}  & 94.1\% & 94.4\% \\
Operations with significant digits & \textbf{80.4\%} & 72.9\% & 67.4\% \\
Equation of motion & 8.5\%  & \textbf{12.8\%} & 12.2\% \\
Kinetic energy & \textbf{73.5\%} & 72.4\% & 71.8\% \\ 
Float or sink? &  98\% & 98.7\% & \textbf{98.8\%} \\ 
Series/Parallel combination of resistance & \textbf{88.9\%} & 32.5\% & 45.7\% \\
\bottomrule
\end{tabular}
\caption{Compare various transformers. Here C2C is char-to-char encoding, W2W is word-to-word encoding, Perfor. stands for performer \cite{choromanski2021rethinking}, and Trans. stands for Transformer \cite{vaswani2017}.}
\label{compare}
\end{table}
\normalsize

\begin{figure}[h]
    \centering
\includegraphics[scale=0.18]{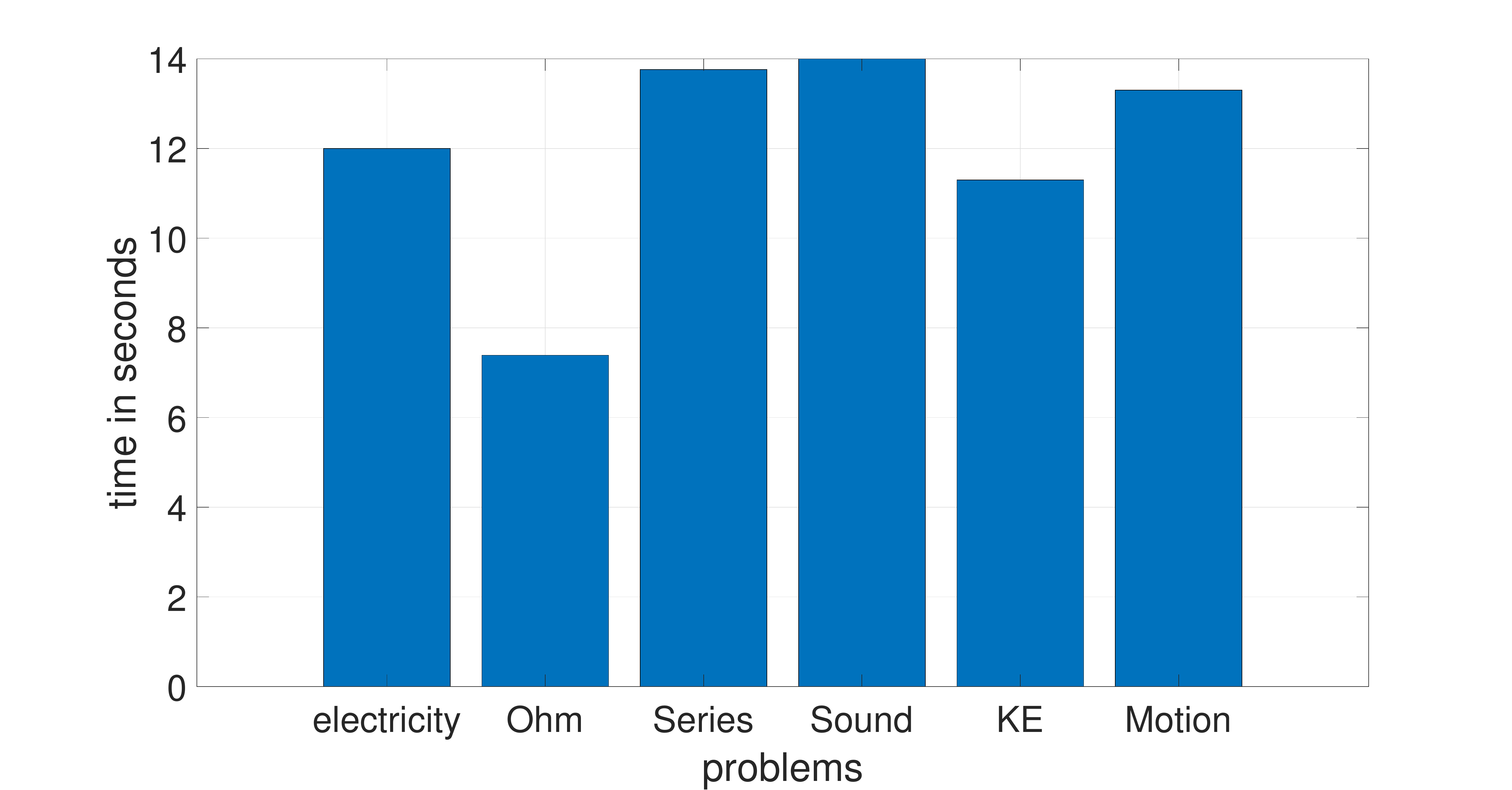}
   \caption{\label{fig:time}Time for generating some datasets. The generator codes are provided in the data repository.}  
\end{figure}

\clearpage

\begin{ack}
This work was done at IIIT, Hyderabad. The authors acknowledge all the support of the institute.
\end{ack}

\end{document}